\documentclass[leqno]{aadmbook}
\usepackage{amsthm}
\usepackage{amsfonts}
\usepackage{amsmath}
\usepackage{graphicx}

\usepackage{url}

\usepackage{bigints} 

\usepackage{arydshln}
\def\nfrac#1#2{\mbox{\small $\displaystyle\frac{#1}{#2}$}}
\def\nnfrac#1#2{\mbox{\footnotesize $\displaystyle\frac{#1}{#2}$}}
\def\nnnfrac#1#2{\mbox{\scriptsize $\displaystyle\frac{#1}{#2}$}}

\begin{document}

\thispagestyle{plain}

$\;$

\noindent
\qquad {\tt   arXiv preprint}
\vspace{2cc}

\noindent
{\small Applicable Analysis and Discrete Mathematics; \\
https://doi.org/10.2298/AADM220526027P}

\vspace*{23.0 mm}


\oddsidemargin 16.5mm
\evensidemargin 16.5mm





\begin{center}

{\large\bf  H\" UGELSCH\" AFFER EGG CURVE AND SURFACE
\rule{0mm}{6mm}\renewcommand{\thefootnote}{}
\footnotetext{\scriptsize ${}^{\ast}$Corresponding author. Maja Petrovi\' c}
\footnotetext{

\rule{2.4mm}{0mm}Keywords and Phrases. egg shaped curve/surface, H\" ugelsch\" affer's model, elliptic integral 
}}

\vspace{1cc}
{\large\it Maja Petrovi\' c\;${}^{\ast}$ and Branko Male\v sevi\' c}

\vspace{1cc}
\parbox{24cc}{{\small

In this paper we consider H\" ugelsch\" affer cubic curves which are generated using appropriate geometric constructions.
The main result of this work is the mode of explicitly calculating the area of the egg-shaped part of the cubic curve using elliptic integrals.
In this paper, we also analyze the H\" ugelsch\" affer surface of cubic curves for which we provide new forms of formulae for the volume and surface area of the egg-shaped part.
Curves and surfaces of ovoid shape have wide applicability in aero-engineering and construction, and are also of biologic importance.
With respect to this, in the final section, we consider some examples of the real applicability of this H\" ugelsch\" affer model.
}}
\end{center}

\bigskip

\bigskip

\begin{center}
{\bf 1. INTRODUCTION}
\end{center}

The starting point is the geometric construction of planar egg-shaped curves introduced by the German engineer {\sc Fritz H\" ugelsch\" affer} in 1944, \cite{XYZ 1944}, \cite{Schmidbauer 1948}, \cite{Schmidbauer 1949}.

\begin{figure}[htb]
\vspace{-0.3cc}
  \centerline{
    \includegraphics*[width=0.6\textwidth]{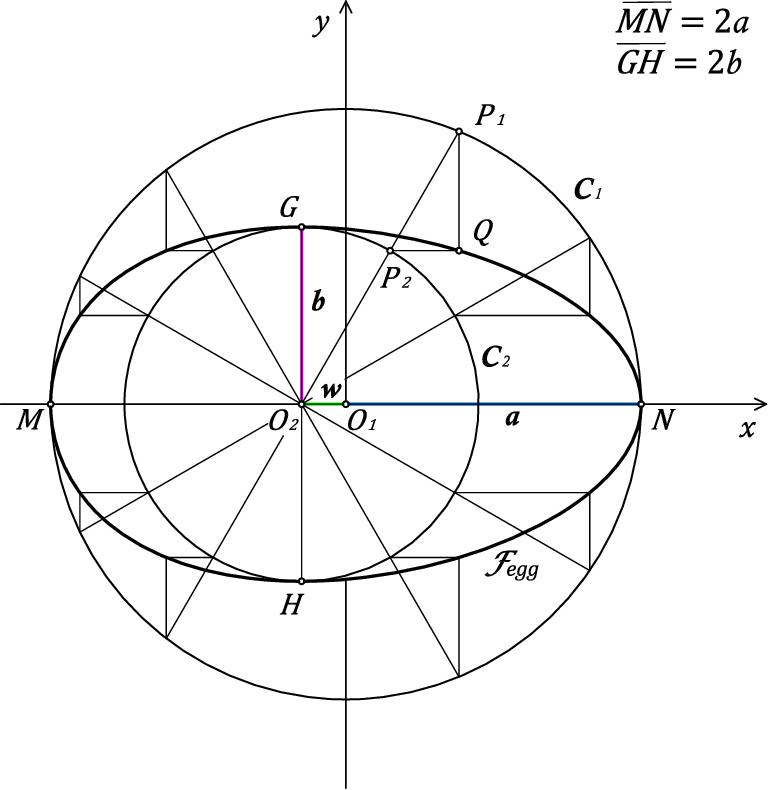}
  }
\caption{H\" ugelsch\" affer's  construction of an egg curve}
\label{Fig.1}       
\end{figure}

\bigskip
{\em H\" ugelsch\" affer curve definition.} Let there exist two circles, in an Euclidean plane:
$$
\mathcal{C}_1: \quad  X^2+Y^2=a^2
$$
and
$$
\mathcal{C}_2: \quad (X+w)^2+Y^2=b^2,
$$
where $a, b, w \ge 0$. For the point $P_1=(X_1,Y_1) \!\in\! \mathcal{C}_1$ we define the point $P_2=(X_2,Y_2) \!\in\! \mathcal{C}_2$ by
$P_2 = O_2P_1 \cap \mathcal{C}_2$, where $O_2=(-w,0)$ is the center of circle $\mathcal{C}_2$, see Figure \ref{Fig.1}.
The locus of points $Q = (x,y)$ where $x=X_1$ and $y=Y_2$, satisfies the cubic algebraic equation
\begin{equation}
\label{F}
\mathcal{F}: \quad 2wxy^2+b^2x^2+(a^2+w^2)y^2-a^2b^2=0,
\end{equation}
as shown in \cite{MP 2010}, see also \cite{Ferreol 2009}.
The H\" ugelsch\" affer curve is defined by the locus of points $Q = (x, y)$. Let us emphasize that the previous definition of the cubic algebraic equation
is consistent with Newton’s hyperbolism \cite{Newton 1860}; for details about hyperbolism see considerations in \cite{Ferreol 2009H}, \cite{MP 2010} and \cite{moNG 2010a}.

\medskip
Let us assume that $a, b, w \ge 0$. We consider the cases where some of the parameters $a, b, w$ equal zero.
If $b=0$, the algebraic curve $\mathcal{F}$ degenerates into an union of lines.
The cases where $a=0$ and/or $w=0$ are of special interest. If $a = 0$ and $w \neq 0$ the algebraic curve $\mathcal{F}$ becomes a mixed cubic \cite{Heinrich 1908}, \cite{MP 2010}, \cite{Ferreol 2001}.
This type of curve was studied by G. de Longchamps in 1890, who named it. If $a \neq 0$ and $w = 0$ the algebraic curve $\mathcal{F}$ becomes an ellipse.
If $a=w=0$ the algebraic curve $\mathcal{F}$ degenerates into an union of lines.
In all other cases ($b \neq 0 \wedge w \neq 0 \wedge a \neq w$) the algebraic curve is non-degenerated algebraic curve of the third order
and in this paper we consider only such curves. Let us note that if $w=a$ the algebraic curve $\mathcal{F}$ degenerates into an union of a line and a parabola, \cite{MP 2010}.

\medskip
Furthermore, let $w>0$. Solving the algebraic equation (\ref{F}) for $y$, two functions are obtained:
\begin{equation}
\label{f1}
y
\!=\!
f_{1,2}(x)
\!=\!
\pm\,
b\mbox{\small $\displaystyle \sqrt{\frac{a^2-x^2}{2wx\!+\!a^2\!+\!w^2}}$}
=
\pm\,
H
\mbox{\small $\displaystyle \sqrt{\frac{(\alpha\!-\!x)(x\!-\!\beta)}{x-\gamma}}$}\!: \left(-\infty,\gamma\right) \cup \left[-a,a\right]\! \longrightarrow R
\end{equation}
where
\begin{equation}
\label{ABC}
\alpha
=
a,\;
\beta
=
-a,\;
\gamma
=
-\frac{a^2+w^2}{2w},\;
H
=
\displaystyle\frac{b}{\sqrt{2w}}.
\end{equation}
as considered in \cite{MP 2010}. It is easy to check that $\gamma <-a \Leftrightarrow (a-w)^2>0$.

\begin{figure}[!h]
  \centerline{
    \includegraphics*[width=\textwidth]{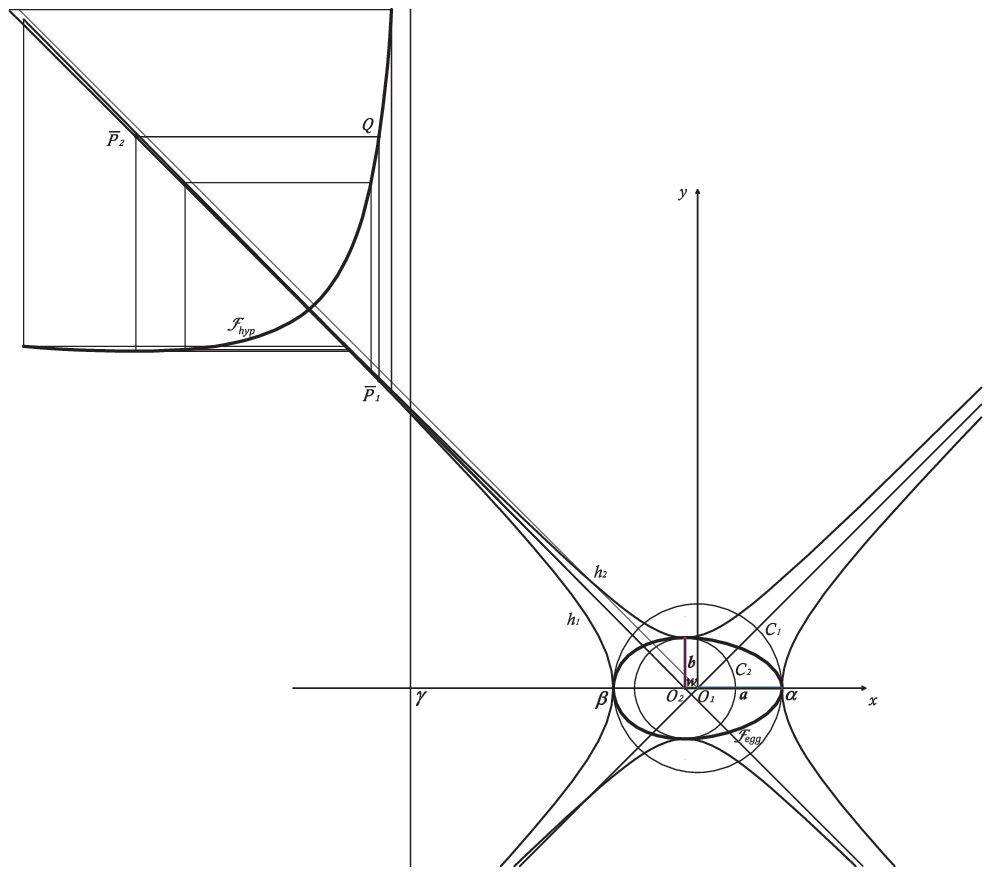}
  }
\caption{H\" ugelsch\" affer's construction of the egg-shaped part and the Petrovi\' c construction of the hyperbolic part of a cubic curve}
\label{Fig.2}       
\end{figure}

Functions $f_{1,2}$ over $[-a,a]$ determine the egg-shaped part of the cubic curve which we denote with $\mathcal{F}_{egg}$.
Functions $f_{1,2}$ over $(-\infty,\gamma)$ determine the hyperbolic part of the cubic curve which we denote with $\mathcal{F}_{hyp}$, see Figure \ref{Fig.2}.
The geometric construction of the curve $\mathcal{F}_{hyp}$, i.e. the extended hyperbolic part of the cubic curve, using two hyperbolae and Newton's hyperbolism
was first given in the magister's thesis \cite{MP 2010} by M. Petrovi\' c. It is noted that $\mathcal{F} = \mathcal{F}_{egg} \cup \mathcal{F}_{hyp}$.

Various analytical properties of the functions $f_{1,2}$ given by (\ref{f1}) are considered in \cite{MP 2010}, \cite{moNG 2010a}.
It is simple to show the following result which is used in the next section. It holds that:

\smallskip\noindent
{\bf Theorem 1.} {\em Let \mbox{\boldmath $u$} be the abscissa of the maximum of function $y = f_1(x)$ over $[-a,a]$. Then$:$
\begin{equation}
\label{Uu}
\mbox{\boldmath $u$}
=
\left\{
\begin{array}{ccc}
-w     &:& w<a, \\[1.0 ex]
-a     &:& w=a, \\[1.0 ex]
-a^2/w &:& w>a.
\end{array}
\right.
\end{equation}}

\smallskip\noindent
{\bf Note 1.} {\sl The case where $w\!=\!a$ leads to the degeneration of the cubic curve ${\cal F}$, so we do not take this case of the previous theorem into further consideration.}

\bigskip

\begin{center}
{\bf 2. MAIN RESULTS}
\end{center}

\bigskip
\noindent
{\bf \boldmath 2.1. H\" ugelsch\" affer's egg-shaped curve}

\bigskip
\noindent
The starting point is the area of the surface bounded by the cubic curve ${\cal F}_{egg}$ which is given by the integral:
\smallskip
$$
\mathcal{A}_{egg}
=
2\displaystyle\int\limits_{\beta}^{\alpha}{f_1(x)\, dx}.
$$
\smallskip
The value $\mathcal{A}_{egg}$ is obtained using the well-known integrals \cite{Ryzhik 2015}:

\smallskip\noindent
{\it The elliptic integral of the first kind}
$$
F(\theta ,p) = \int\limits_{0}^{\sin(\theta)}  \frac{dt}{\sqrt{(1-t^2)(1-p^2t^2)}}, \;\;\;\; 0 \leq p^2 \leq 1 \;\;\;\; \wedge \;\;\;\; 0 \leq \sin(\theta) \leq 1,
$$

\smallskip\noindent
and {\it the elliptic integral of the second kind}
$$
E(\theta ,p) = \int\limits_{0}^{\sin(\theta)}  \frac{\sqrt{1-p^2t^2}}{{\sqrt{1-t^2}}}\,\,  dt, \;\;\;\; 0 \leq p^2 \leq 1 \;\;\;\; \wedge \;\;\;\; 0 \leq \sin(\theta) \leq 1.
$$
Furthermore, let there be functions
$$
\varphi(x) = \sqrt{\frac{(\alpha-\gamma)(x-\beta)}{(\alpha-\beta)(x-\gamma)}}
\;\;\;\; \wedge \;\;\;\;
\psi(x) = \sqrt{\frac{\alpha-x}{\alpha-\beta}},
$$
\smallskip\noindent
where the parameters $\alpha, \beta$ and $\gamma$ are determined with (\ref{ABC}) and where $x \in [\beta,\alpha]$.
In the next theorem, we will use the following lemma.

\bigskip\noindent
{\bf Lemma 1.}
{\em
$$
\varphi(x)
=
\psi(x)
\quad\Longleftrightarrow\quad
x = - w
\;\vee\;
x = - \frac{a^2}{w}.
$$}
\smallskip\noindent
The following statement is true:

\smallskip\noindent
{\bf Theorem 2.}
{\em
It holds that
\begin{equation}
\label{Aegg}
\mathcal{A}_{egg}
=
\mathcal{A}_1
+
\mathcal{A}_2
\end{equation}
where
$$
\begin{array}{rcl}
\mathcal{A}_1
\!\!&\!\!=\!\!&\!\!
\displaystyle\frac{4}{3} H \sqrt{\alpha-\gamma} \, {\Big (} (\alpha+\beta-2\gamma)\cdot E(\kappa ,p) - 2(\beta-\gamma)\cdot F(\kappa ,p) {\Big )} +                 \\[2.5 ex]
\!\!&\!\! \!\!&\!\!
\displaystyle\frac{4}{3} H (\mbox{\boldmath $u$}+\gamma-\alpha-\beta)\sqrt{\frac{(\alpha-\mbox{\boldmath $u$})(\mbox{\boldmath $u$}-\beta)}{\mbox{\boldmath $u$}-\gamma}}
\end{array}
$$
and
$$
\begin{array}{rcl}
\mathcal{A}_2
\!\!&\!\!=\!\!&\!\!
\displaystyle\frac{4}{3} H \sqrt{\alpha-\gamma} \, {\Big (} (\alpha+\beta-2\gamma)\cdot E(\lambda ,p) - 2(\beta-\gamma)\cdot F(\lambda ,p) {\Big )} -               \\[2.5 ex]
\!\!&\!\! \!\!&\!\!
\displaystyle\frac{4}{3} H \sqrt{(\alpha-\mbox{\boldmath $u$})(\mbox{\boldmath $u$}-\beta)(\mbox{\boldmath $u$}-\gamma)}
\end{array}
$$
with the parameters$:$
$$
p = \sqrt{\frac{\alpha-\beta}{\alpha-\gamma}}, \quad
\kappa = \arcsin \varphi(\mbox{\boldmath $u$}), \quad
\lambda = \arcsin \psi(\mbox{\boldmath $u$})
$$
and where \mbox{\boldmath $u$} is the abscissa of the maximum of function $y = f_1(x)$  $($see eq.~$(\ref{Uu}))$.}

\begin{figure}[!h]
\vspace*{-4.0 mm}
  \centerline{
    \includegraphics*[width=0.55\textwidth]{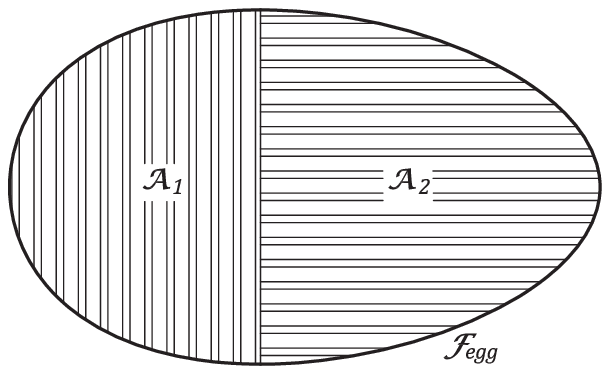}
  }
\caption{Areas $\mathcal{A}_1$ and $\mathcal{A}_2$ of the H\" ugelsch\" affer egg curve}
\label{Fig.3}       
\end{figure}

\bigskip
\noindent
{\bf Proof.} Note that $\alpha > \beta > \gamma$ and it follows that for $p$ it holds
$$
0 < p^2 =
\frac{\alpha-\beta}{\alpha-\gamma} = 1 - \frac{\beta-\gamma}{\alpha-\gamma} < 1.
$$
 Starting from {\boldmath $u$} as the abscissa of the maximum of the function for parameters
$
\kappa = \arcsin \varphi(\mbox{\boldmath $u$})
$
and
$\lambda = \arcsin \psi(\mbox{\boldmath $u$})
$
it holds that
$$
\kappa = \lambda = \arcsin \sqrt{\frac{\alpha-\mbox{\boldmath $u$}}{\alpha-\beta}},
$$
in accordance with Theorem 1. and Lemma 1. Additionally, the expression
$$
0
<
\sin \lambda
=
\sqrt{\frac{\alpha-\mbox{\boldmath $u$}}{\alpha-\beta}}
=
\sqrt{1-\frac{\mbox{\boldmath $u$}-\beta}{\alpha-\beta}}
<
1,
$$
is satisfied, because $\alpha > \mbox{\boldmath $u$} > \beta > \gamma$. Because of this, the inequality
$$
0
<
\sin \kappa
<
1.
$$
is also correct. First, let us consider
$$
I_1
=
\displaystyle\int\limits_{\beta}^{\mbox{\footnotesize \boldmath $u$}}{f_1(t)\, dt}
=
H \displaystyle\int\limits_{\beta}^{\mbox{\footnotesize \boldmath $u$}}{\sqrt{\frac{(\alpha-t)(t-\beta)}{t-\gamma}}\, dt.}
$$
Then, from \cite{Ryzhik 2015} according to formula 3.141/34 it holds that
$$
\begin{array}{rcl}
I_1
\!\!&\!\!=\!\!&\!\!
\displaystyle\frac{2}{3} H \sqrt{\alpha-\gamma} \, {\Big (} (\alpha+\beta-2\gamma)\cdot E(\kappa ,p) - 2(\beta-\gamma)\cdot F(\kappa ,p) {\Big )} +                 \\[2.5 ex]
\!\!&\!\! \!\!&\!\!
\displaystyle\frac{2}{3} H (\mbox{\boldmath $u$}+\gamma-\alpha-\beta)\sqrt{\frac{(\alpha-\mbox{\boldmath $u$})(\mbox{\boldmath $u$}-\beta)}{\mbox{\boldmath $u$}-\gamma}}
\end{array}
$$
because the condition $\alpha \geq \mbox{\boldmath $u$} > \beta > \gamma$ is satisfied.
Second, let us consider
$$
I_2
=
\displaystyle\int\limits_{\mbox{\footnotesize \boldmath $u$}}^{\alpha}{f_1(t)\, dt}
=
H \displaystyle\int\limits_{\mbox{\footnotesize \boldmath $u$}}^{\alpha}{\sqrt{\frac{(\alpha-t)(t-\beta)}{t-\gamma}}\, dt}.
$$
Then, from \cite{Ryzhik 2015} according to formula 3.141/35 it holds that
$$
\begin{array}{rcl}
I_2
\!\!&\!\!=\!\!&\!\!
\displaystyle\frac{2}{3} H \sqrt{\alpha-\gamma} \, {\Big (} (\alpha+\beta-2\gamma)\cdot E(\lambda ,p) - 2(\beta-\gamma)\cdot F(\lambda ,p) {\Big )} -               \\[2.5 ex]
\!\!&\!\! \!\!&\!\!
\displaystyle\frac{2}{3} H \sqrt{(\alpha-\mbox{\boldmath $u$})(\mbox{\boldmath $u$}-\beta)(\mbox{\boldmath $u$}-\gamma)}
\end{array}
$$
because the condition $\alpha > \mbox{\boldmath $u$} \geq \beta > \gamma$ is satisfied.
Therefore, the claim of the theorem follows from
$$
\mathcal{A}_1 = 2 I_1 \;\wedge\; \mathcal{A}_2 = 2 I_2 \quad \Longrightarrow \quad \mathcal{A}_{egg} = \mathcal{A}_1 + \mathcal{A}_2.$$
Also see Figure \ref{Fig.3}.

\vspace*{-12.0 mm}
$\,$
$$
 \eqno \Box
$$

\bigskip
\noindent
Let us emphasize that if $\mbox{\boldmath $u$}=-w$, then
\begin{equation}
\label{Case_1_S_1}
\begin{array}{rcl}
\mathcal{A}_{1}
\!\!&\!=\!\!&\!\!
\displaystyle\frac{2(a+w)b}{3w}
{\Bigg (}
\displaystyle\frac{a^2+w^2}{w}
E{\bigg (}\!\!\arcsin \sqrt{\displaystyle\frac{a+w}{2a}},\displaystyle\frac{2\sqrt{aw}}{a+w}\,{\bigg )} -                \\[2.5 ex]
\!\!&\! \!\!&\!\!
\displaystyle\frac{(a-w)^2}{w}
F{\bigg (}\!\!\arcsin \sqrt{\displaystyle\frac{a+w}{2a}},\displaystyle\frac{2\sqrt{aw}}{a+w}\,{\bigg )} \!
{\Bigg )}
-
\displaystyle\frac{2b}{3w}(a^2+3w^2)
\end{array}
\end{equation}
and
\begin{equation}
\label{Case_1_S_2}
\begin{array}{rcl}
\mathcal{A}_{2}
\!\!&\!=\!\!&\!\!
\displaystyle\frac{2(a+w)b}{3w}
{\Bigg (}
\displaystyle\frac{a^2+w^2}{w}
E{\bigg (}\!\!\arcsin \sqrt{\displaystyle\frac{a+w}{2a}},\displaystyle\frac{2\sqrt{aw}}{a+w}\,{\bigg )} -                \\[2.5 ex]
\!\!&\! \!\!&\!\!
\displaystyle\frac{(a-w)^2}{w}
F{\bigg (}\!\!\arcsin \sqrt{\displaystyle\frac{a+w}{2a}},\displaystyle\frac{2\sqrt{aw}}{a+w}\,{\bigg )} \!
{\Bigg )}
-
\displaystyle\frac{2b}{3w}(a^2-w^2).
\end{array}
\end{equation}
Also, let us emphasize that if $\mbox{\boldmath $u$}=-\displaystyle\frac{a^2}{w}$, then
\begin{equation}
\label{Case_2_S_1}
\begin{array}{rcl}
\mathcal{A}_{1}
\!\!&\!=\!\!&\!\!
\displaystyle\frac{2(a+w)b}{3w}
{\Bigg (}
\displaystyle\frac{a^2+w^2}{w}
E{\bigg (}\!\!\arcsin \sqrt{\displaystyle\frac{a+w}{2w}},\displaystyle\frac{2\sqrt{aw}}{a+w}\,{\bigg )} -                \\[2.5 ex]
\!\!&\! \!\!&\!\!
\displaystyle\frac{(a-w)^2}{w}
F{\bigg (}\!\!\arcsin \sqrt{\displaystyle\frac{a+w}{2w}},\displaystyle\frac{2\sqrt{aw}}{a+w}\,{\bigg )} \!
{\Bigg )}
-
\displaystyle\frac{2ab}{3w^2}(3a^2+w^2)
\end{array}
\end{equation}
and
\begin{equation}
\label{Case_2_S_2}
\begin{array}{rcl}
\mathcal{A}_{2}
\!\!&\!=\!\!&\!\!
\displaystyle\frac{2(a+w)b}{3w}
{\Bigg (}
\displaystyle\frac{a^2+w^2}{w}
E{\bigg (}\!\!\arcsin \sqrt{\displaystyle\frac{a+w}{2w}},\displaystyle\frac{2\sqrt{aw}}{a+w}\,{\bigg )} -                \\[2.5 ex]
\!\!&\! \!\!&\!\!
\displaystyle\frac{(a-w)^2}{w}
F{\bigg (}\!\!\arcsin \sqrt{\displaystyle\frac{a+w}{2w}},\displaystyle\frac{2\sqrt{aw}}{a+w}\,{\bigg )} \!
{\Bigg )}
-
\displaystyle\frac{2ab}{3w^2}(-a^2+w^2).
\end{array}
\end{equation}
Notice that from (\ref{Case_1_S_1}) and (\ref{Case_1_S_2}) it follows
$$
\mathcal{A}_1 + \mathcal{A}_2 = \mathcal{A}_{egg} \quad\mbox{and}\quad \mathcal{A}_2 - \mathcal{A}_1 = \frac{8}{3}bw;
$$
and that from (\ref{Case_2_S_1}) and (\ref{Case_2_S_2}) it follows
$$
\mathcal{A}_1 + \mathcal{A}_2 = \mathcal{A}_{egg} \quad\mbox{and}\quad \mathcal{A}_2 - \mathcal{A}_1 = \frac{8}{3}\frac{a^3b}{w^2}.
$$


\newpage
\noindent
{\bf \boldmath 2.2. H\" ugelsch\" affer's surface}

\bigskip
\noindent
{\em $1.$ The implicit form of the {\sc H\" ugelsch\" affer}'s surface.}
The 3D surface generated by rotating the curve ${\cal F}$ given by (\ref{F}) around the $x$-axis for $360^{\circ}$ has
the following implicit form
\begin{equation}
\label{egg22}
\frac{x^2}{a^2}+\frac{y^2 + z^2}{b^2}g(x)=1, \qquad g(x)=1+\frac{2wx+w^2}{a^2}.
\end{equation}

\begin{figure}[!h]
  \centerline{
    \includegraphics*[width=\textwidth]{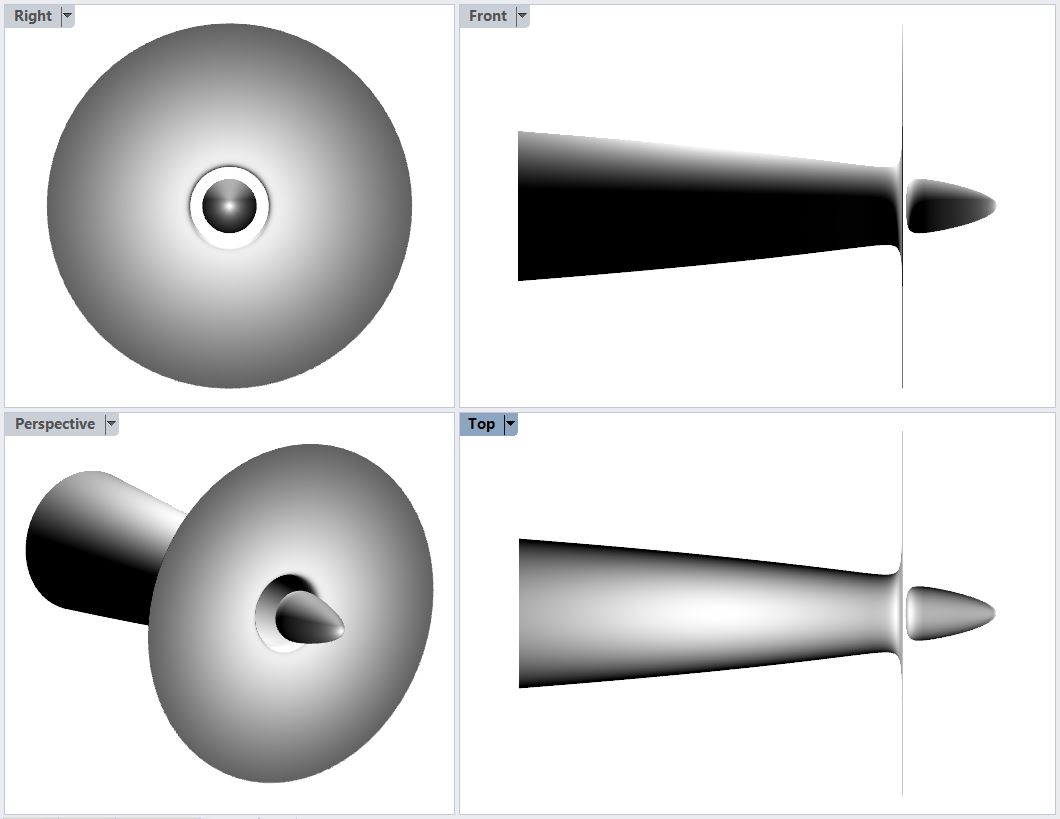}
  }
\caption{A H\" ugelsch\" affer surface (top, front, right and axonometric view)}
\label{Fig.4}       
\end{figure}

\noindent
In Figure \ref{Fig.4}, a 3D surface generated using the cubic curve $\mathcal{F}$ with appropriate parameters $a$, $b$ and $w$ is shown.

\bigskip
\noindent
{\em $2.$ The explicit form of the {\sc H\" ugelsch\" affer} surface.}
Solving the equation (\ref{egg22})~for~$z$, the explicit form of the {\sc H\" ugelsch\" affer} surface is obtained
$$
z=z(x,y)=\pm \sqrt{\frac{b^2(a^2 - x^2)}{2wx+w^2+a^2}-y^2},
$$
with domain
$
D_{xy} = {\big \{}\,\! (x,y) \!\in\! R^2 \;|\; f_1(x) \leq y \leq f_2(x)  : x \!\in\! (-\infty, \gamma) \!\cup\! [-a,a] \,\! {\big \}}.
$

\break

\noindent
{\em $3.$ The volume of the surface generated by rotating the curve $\mathcal{F}_{egg}$ around the $x$-axis.}
It is simple to show that the volume of the surface generated by rotating the curve $\mathcal{F}_{egg}$ around the $x$-axis is given by
\begin{equation}
\label{Volume}
V_{egg}
=
\pi \!\displaystyle\int\limits_{\beta}^{\alpha}{\!f_1^2(x)\, dx}
=
\frac {\pi b^2}{4w^3} \left( (a^2-w^2)^2 \ln \left |\frac{a-w}{a+w}\right |+ 2aw (a^2+w^2) \right).
\end{equation}

\noindent
\quad For the previous expression for surface volume it holds:

\medskip
\noindent
$(i)$
if observing $w$ in the following limit process, the volume of a spheroid is obtained
$$
V_{spheroid}
=
\lim\limits_{w \rightarrow 0+}{V_{egg}}
=
\frac{4}{3} ab^2\pi\,.
$$
Note that, from this, it is obtained, for the value $b=a$:
$$
V_{sphere}
=\frac{4}{3} a^3\pi\,,
$$
which represents the formula for the volume of a sphere of radius $a$.

\medskip
\noindent
$(ii)$
if $w=a$ and $a \neq 0$, then the volume of a paraboloid is obtained:
$$
V_{paraboloid}
=ab^2\pi.
$$
{\em $4.$ The surface area of the surface generated by rotating the curve $\mathcal{F}_{egg}$ around the $x$-axis for $360^{\circ}$.}
 The egg-shaped surface has the following formula for calculating the area of the 3D surface:

$$
{\cal S}_{egg}
=
2\pi \!\displaystyle\int\limits_{\beta}^{\alpha} {\!f_1(x) \sqrt{1+\left( f_1^{\prime}(x)\right) ^2} \,\, dx}.
$$
Substituting $x = t + \gamma$ it is obtained:
$$
{\cal S}_{egg}
=
\nnfrac{b \pi}{4 w^2} \!\!\!\displaystyle\int\limits_{-a-\gamma}^{a-\gamma} \!\!\!\nnfrac{\sqrt{Q_5(t)}}{t^2}\,\, dt \,,
$$
with the polynomial of the fifth degree
$$
Q_5(t)
=
a_5 t^5 + a_4 t^4 + a_3 t^3 + a_2 t^2 + a_1 t + a_0
$$
and the coefficients
$$
\begin{array}{rcl}
a_5 \!\!&\!\!=\!\!&\!\! -32 w^3                                    \\[1.5 ex]
a_4 \!\!&\!\!=\!\!&\!\! 4 w^2(8a^2+b^2+8w^2)                       \\[1.5 ex]
a_3 \!\!&\!\!=\!\!&\!\! -8 w (a^2-w^2)^2                           \\[1.5 ex]
a_2 \!\!&\!\!=\!\!&\!\! -2 b^2 (a^2-w^2)^2                         \\[1.5 ex]
a_1 \!\!&\!\!=\!\!&\!\! 0                                          \\[1.5 ex]
a_0 \!\!&\!\!=\!\!&\!\! \nnfrac{b^2}{4w^2} (a^2-w^2)^4.
\end{array}
$$
Based on the previous, it follows that
$$
\label{zvezdica}
\!\!\!\!\!\!\!\!
\begin{array}{l}
(12) \;\; {\cal S}_{egg} =                                                    \\[2.0 ex]
\!\!\! = \nnnfrac{b \pi}{4 w^2} \!\!\!\!\!\!\!\!
\mathop{\mbox{$\displaystyle \bigintss$}}\limits_{\frac{(a-w)^2}{2w}}^{\frac{(a+w)^2}{2w}}{
\!\!\!\!\!\!\!\nnnfrac{\sqrt{
-32 w^3 t^5
\!+\!
4 w^2 (8a^2\!+\!b^2\!+\!8w^2) t^4
\!-\!
8 w (a^2\!-\!w^2)^2 t^3
\!-\!
2 b^2 (a^2\!-\!w^2)^2 t^2
\!+\!
\nnnfrac{b^2}{4 w^2}(a^2\!-\!w^2)^4
}}{t^2}
dt.
}
\end{array}
$$ \setcounter{equation}{12}
\noindent
Let us consider some approximative methods for calculating the integral (12). Utilizing the Simpson's quadrature rule for $n=3$ points:
$$
{\cal S}_{egg}
\approx
{\cal S}_{{egg}_{n=3}}
=
\nnfrac{a}{3}\nnfrac{b \pi}{4 w^2}
\left(\nnfrac{\sqrt{Q_5(-a-\gamma)}}{(-a-\gamma)^2}+4\nnfrac{\sqrt{Q_5(-\gamma)}}{\gamma^2}+\nnfrac{\sqrt{Q_5(a-\gamma)}}{(a-\gamma)^2}\right)=
$$
$$=
\nnfrac{\pi  a}{3} 
\left(\nnfrac{2 ab^2}{(a-w)^2}+\nnfrac{8 ab}{(a^2+w^2)^2}\sqrt{(a^2+w^2)^3+a^2b^2w^2}+\nnfrac{2 ab^2}{(a+w)^2}\right)
$$
an approximative formula for the area of the {\sc H\" ugelsch\" affer} egg-shaped surface is obtained.

\medskip
\noindent
\quad For the previous expression for the area of the surface, it holds that:

\medskip
\noindent
$(i)$
if observing $a \neq 0$ and $w \rightarrow 0+$ an approximative formula for the area of a spheroid is obtained:
$$
{\cal S}_{spheroid}
\approx
\nnfrac{4b \pi}{3}(b+2a)\,.
$$
This approximative formula is also given in \cite{Spheroids 1960}. Note that for the value $b=a$ it holds that
$$
{\cal S}_{sphere}
=
4a^2 \pi \,,
$$
which represents the formula for the area of a sphere of radius $a$.

\medskip
\noindent
$(ii)$
if observing $w=a$ and $a \neq 0$, from $(12)$ a well-known formula for the area of a paraboloid is obtained:
$$
{\cal S}_{paraboloid}
=
{
\nfrac{b \pi}{2 a^2} \!\!\!\!
\mathop{\mbox{$\displaystyle \bigintsss$}}\limits_{\!\!\!\!\!\!\!0}^{\;\;\;\;2a}{
\!\!\!\!\nfrac{\sqrt{
-8a^3 t^{5}
+
a^2(8a^2\!+\!b^2\!+\!8a^2) t^4
}}{t^2} \, dt}
}
=
\nfrac{b \pi}{24a^2}(\left(16a^2+b^2\right)^{\!3/2}-b^3) \, .
$$

\newpage
\begin{center}
{\bf 3. APPLICATIONS}
\end{center}

The construction of the {\sc H\" ugelsch\" affer} egg-shaped curve model was established in the mid-20th century for the purposes of aero-engineering \cite{XYZ 1944}, \cite{Schmidbauer 1948}, \cite{Schmidbauer 1949}. Even today, this model is used for the fuselage cross-sections of gliders, \cite{Boermans 2004}:
"The fuselage cross-sections of the ASW27 are described by a so-called 'H\" ugelsch\" affer-Egg-Curve' which has the special feature of a continuous curvature
on its circumference, a prerequisite for a smooth pressure distribution and undisturbed boundary layer development" and according to \cite{Boermans 2006}:
"For an undisturbed boundary layer development, continuity of curvature is required in flow direction. This is guaranteed by deriving the top, bottom, and line of largest width from airfoil shapes, and using H\" ugelsch\" affer curves (deformed ellipses) for the fuselage cross-sections".

In the papers \cite{Majamath 2019}, \cite{EggSurf 2015}, \cite{moNG 2010b}, \cite{moNG 2008}, \cite{ICEGD 2013} and \cite{ICEGD 2011}
some mathematical properties of the {\sc H\" ugelsch\" affer} curve are considered, as well as some of its applications in architecture.

The applications of egg-shaped curves are also significant in civil engineering, and especially in hydrology, \cite{Kolkata 2013}, \cite{FlowMI 2011},  \cite{WaterS 2021}, \cite{Water 2016}, \cite{Water 2019}, \cite{WaterSIE 2017},  \cite{WaterST 2013}, \cite{WaterER 2022}. Sewer Systems with Egg-Shaped Pipe Cross Sections are gaining traction in cities because of the availability of modern technologies. Egg-shaped or ovoid cross-sections are also common in old brick-built sewerage systems.
Several cities around the world including London, Paris, Hamburg, New York, Los Angeles, Cleveland, Minneapolis, Newark, Kolkata, Mumbai and Delhi have such sewer systems.
The egg-shaped cross section of these pipes is achieved by the composition of circular arcs.
Note that in India, according to a standard from 1976, sewer pipes with an egg-shaped trifocal curve profile (egglipse section) are used, \cite{Indian 1976}.

In the next part of this section, we state some numerical examples.

\smallskip
\textbf{Numerical example 1}:
For the city of Cleveland, Ohio, USA, in \cite{Cleveland 2008}, a graphical model of egg-shaped curves for sewer channels with length $L=H$, shape index  $B/L=\frac{2+\sqrt{3}}{3+\sqrt{3}}$ and offset from the bottom to the point where the channel is widest $S=\frac{5}{2} \frac{L}{3+\sqrt{3}}$ is given. The profile of the sewer pipe is an egg-shaped curve constructed through composition of circular arcs. Based on this model, a {\sc H\" ugelsch\" affer} egg-shaped curve model can be formed with parameters $L$, $B$ and $w=\frac{2-\sqrt{3}}{3+\sqrt{3}} \frac{L}{2}$. In the following table, a comparison of values for area $A$ according to \cite{Cleveland 2008} (see Table: Egg-Shaped Sewers; File No.$\,$73M) and values for the area of the corresponding {\sc H\" ugelsch\" affer} curve $\mathcal{A}$ calculated according to our paper (see Table 1: column 7) is given. Let us emphasize that based on this data, the coefficient of determination is $R^2 = 0.999179\ldots$  The value of the coefficient of determination of $R^2 \approx 1$ confirms that this {\sc H\" ugelsch\" affer} egg-shaped curve model can be a valid approximation of the egg-shaped curve considered in \cite{Cleveland 2008}.

Analogous to the previous example, a {\sc H\" ugelsch\" affer} egg-shaped curve model can be generated according to the profiles of egg-shaped pipe sewers obtained through composition of circular arcs with the shape index $B/L=2/3$ or $B/L=3/4$ as given in \cite{Tables 2006}.

{$$
\renewcommand{\arraystretch}{1.38}
\begin{array}{|c||c|c|c|c||c|c|} \hline
No. & L(ft) & B(ft) & S(ft) &  A(ft^2)&  w(ft) & {\bf \mathcal{A}}(ft^2)  \\ \hline \hline
 2 & 2.25 & 1.94 &  1.28 &  3.41 & 0.155 &  3.420108855 \\ \hdashline
 3 & 2.75 & 2.23 &  1.64 &  4.75 & 0.265 &  4.793986508 \\ \hdashline
 4 & 3.23 & 2.55 &  1.71 &  6.37 & 0.095 &  6.46613264 \\ \hdashline
 5 & 3.74 & 2.95 &  1.98 &  8.54 & 0.110 &  8.66154837 \\ \hdashline
 6 & 4.23 & 3.34 &  2.24 & 10.93 & 0.125 & 11.09141523 \\ \hdashline
 7 & 4.69 & 3.70 &  2.48 & 13.43 & 0.135 & 13.62336584 \\ \hdashline
 8 & 5.12 & 4.04 &  2.71 & 16.01 & 0.150 & 16.23882865 \\ \hdashline
 9 & 5.54 & 4.38 &  2.93 & 18.74 & 0.160 & 19.04989185 \\ \hdashline
10 & 5.94 & 4.69 &  3.14 & 21.55 & 0.170 & 21.87112886 \\ \hdashline
11 & 6.33 & 5.00 &  3.35 & 24.47 & 0.185 & 24.84723069 \\ \hdashline
12 & 6.71 & 5.30 &  3.55 & 27.50 & 0.195 & 27.9193156 \\ \hdashline
13 & 7.08 & 5.59 &  3.75 & 30.61 & 0.210 & 31.0701810 \\ \hdashline
14 & 7.44 & 5.88 &  3.94 & 33.80 & 0.220 & 34.3439429 \\ \hdashline
15 & 7.79 & 6.15 &  4.12 & 37.06 & 0.225 & 37.6115462 \\ \hdashline
16 & 8.13 & 6.42 &  4.30 & 40.37 & 0.235 & 40.9764108 \\ \hdashline
17 & 8.47 & 6.69 &  4.48 & 43.81 & 0.245 & 44.4854110 \\ \hdashline
18 & 8.79 & 6.94 &  4.65 & 47.19 & 0.255 & 47.8911598 \\ \hdashline
19 & 9.12 & 7.20 &  4.83 & 50.79 & 0.270 & 51.5497738 \\ \hdashline
20 & 9.43 & 7.45 &  4.55 & 54.31 & 0.275 & 55.1534978 \\ \hline
\end{array}
$$}
\begin{center}
\textit{Table 1}: Cross section areas of Egg-Shaped Sewers with respective parameters
\end{center}

\vspace*{-3.5 mm} $\,$ \hfill $\Box$

\bigskip
In recent years, the application of {\sc H\" ugelsch\" affer} models, as well as other mathematical models for describing the shape of eggs - ovoid forms (curves and surfaces) in the poultry industry and food engineering is on the rise, as evident from papers \cite{BEggC 2018}, \cite{Heck 2008}, \cite{UTEP 2022}, \cite{VMCh 2021}, \cite{Milojevic 2020}, \cite{EggCurve 2017} -- \cite{Narushin 2022b}.

The reason for the use of the {\sc H\" ugelsch\" affer} egg-shaped curve model is its presentation as a non-destructive oomorphological model when calculating the area of a planar curve $\mathcal{A}$ and the volume $V$ and the surface area $S$ of a 3D egg surface \cite{Narushin 2020b},  \cite{Narushin 2021d}, \cite{Narushin 2022}. Note that "oomorphological models" is used as a substitute for "egg-shaped models" \cite{Kostin 1977}.
The authors Narushin et al. used 2D digital images of egg contours to obtain concrete values for the parameters $L=2a$, $B=2b$ and $w$, \cite{Narushin 2020a}, \cite{Narushin 2020b}. According to \cite{Narushin 2020b} the typical length of a hen egg varies between $5$ and $7$~cm, and the ratio of parameters $B$ and $L$ which is called the shape index $SI=B/L$ varies between $0.70$ and $0.78$.
Based on this data and digital imaging of hen eggs Narushin et al. reached the conclusion that the parameter $w$ has a minimum value of $0.021$ and a maximum value of $0.249$ (\cite{Narushin 2020b}, see Table~1). Note that in paper \cite{Milojevic 2020} and \cite{Severa 2013}, the authors consider other ranges for the values of parameters $B$ and $L$.

\medskip
The formulae from the second section of this paper produce correct results for the area of a planar curve $\mathcal{A}$ -- formula (\ref{Aegg}),
the volume $V\!$ -- formula (\ref{Volume}) and the surface area~$\cal S$ of a 3D egg surface -- formula $(12)$. In the following several numerical examples,
a comparison of these formulae with approximative formulae from works
\cite{Narushin 2020a}, \cite{Narushin 2020b} and \cite{Narushin 2021a} is given.

\medskip
\textbf{Numerical example 2}:
In the paper \cite{Narushin 2021a}, the case where $L=2a=6\,cm$, $B/L=0.775$ and $w/L=0.125$ was considered.
In this case, the parameters $a$, $b$ and $w$ for the {\sc H\" ugelsch\" affer} egg-shaped curve model have the following values: $a=L/2=3\,cm$, $b=B/2=4.65/2=2.325\,cm$ and $w=L/8=0.75\,cm$.
Based on the approximative formula given by the authors of \cite{Narushin 2021a}: $\mathcal{A}rea=0.118B^2+0.637LB+0.014L^2$ the area of the plane curve $\mathcal{A}rea=20.827755\,cm^2$ was obtained.
The same authors, in paper \cite{Narushin 2020b}, utilizing Simpson's rule with three points, obtained an approximate value for the area of the plane curve $\mathcal{A}_{Simpson}=21.062316\,cm^2.$
Comparing these values with the values obtained using formulae (\ref{Case_1_S_1}) and (\ref{Case_1_S_2}) as exact results:
$$\mathcal{A}_1 = 8.545026 ...\,cm^2, \,\, \mathcal{A}_2 = 13.195026 ... \,cm^2
 \,\, \Rightarrow \mathcal{A} = \mathcal{A}_1 + \mathcal{A}_2 = 21.740052 ...\, cm^2,$$
we can conclude that
$$
| \mathcal{A} - \mathcal{A}rea | = 0.912297...\,cm^2 \,\,\,\,\,\,\,\,\textrm{and}\,\,\,\,\,\,\,\, | \mathcal{A} - \mathcal{A}_{Simpson} | = 0.677736...\,cm^2. \,\,\,\,\,\,\,\,\,\,\,\,\,\,\,\,\,\,\,\,\,\,\,\Box
$$
\vspace*{-5.5 mm}

\medskip
\textbf{Numerical example 3}: Based on data from Table 2 obtained by the authors of work \cite{Narushin 2020a} through digitally measuring 40 chicken eggs, the mean values of the relevant parameters are: $L=2a=405.81\,pixels,\, B=312.65\,pixels$ and $A=98984.1\,pixels$ i.e. $A=19.05\,cm^2$. Note that in \cite{Narushin 2020a} the following graphical conversion was used: $72.09\,pixels$ in $1\,cm$ and $5197.03\,pixels$ in $1\,cm^2$. Furthermore, from formula $\mathcal{A}rea=0.118B^2+0.637LB+0.014L^2$
of the aforementioned work, the value of the area of the plane curve $\mathcal{A}rea=94660.37552\,pixels$ was obtained, i.e. $\mathcal{A}rea=18.214\,cm^2$.

Observing the equation
$
\mathcal{A} = A
$
it can be concluded that the respective {\sc H\" ugelsch\" affer} egg-shaped curve model has an unique value for the parameter $w$ given by
$w_0=46.678275...\,pixels$ i.e. $w_0=0.6475...\,cm$.
In the work \cite{Narushin 2020b} for the observed data, a mean value of parameter $w=0.12\, cm$ is given. Based on the fact that
$|w_0 - w| > 0.5 \, cm$
it can be concluded that, for these values of parameters $w$ and $w_0$, the contours of the respective egg-shaped curves, despite the area in the {\sc H\" ugelsch\" affer} models being the same, are, geometrically, significantly different.

\hfill $\Box$

\medskip
\textbf{Numerical example 4}: Based on data from Table 1 which authors of the work \cite{Sedghi 2022}, Sedghi and  Ghaderi, obtained by digitally measuring 177 hen eggs at 60 weeks of age, their average values of the relevant parameters are: $L=2a=5.708\,cm, B=4.431\,cm$ and $V=57.458\,cm^3$ (\cite{Sedghi 2022}, see Table 1). Using the following values in our egg-shaped surface model:
 $a=L/2=2.854\,cm$, $b=B/2=2.2155\,cm$ and $V_{egg}=V(=57.458\,cm^3)$, then, according to formula (\ref{Volume}) for the volume of the {\sc H\" ugelsch\" affer} egg-shaped surface model, an unique positive value of parameter $w=0.9138298\ldots\,cm$ can be determined, for which $V_{egg}=V$. For these values of parameters $a$, $b$ and $w$ the surface area of the {\sc H\" ugelsch\" affer} model is $\mathcal{S}_{egg}=73.61192\,cm^2$ which corresponds to some predicted means for $S$ (\cite{Sedghi 2022}, see Table 3, Code A-K). \hfill $\Box$

\medskip
The numerical values in this paper were determined using the computer algebra system Maple. For the purposes of this paper, a web applet dedicated to all the numerical values considered here was also developed \cite{Applet 2022}.

\begin{center}
{\bf 4. CONCLUSION}
\end{center}

Egg-shaped curves and surfaces have wide applications in engineering (air, civil, food).
For the calculation of $A$, $V$ and $S$ of egg-shaped curves and surfaces, mostly numerical approximations were used, such as the ones presented in works \cite{VMCh 2021}, \cite{Narushin 2020a}, \cite{Narushin 2021b}, \cite{Paganelli 1974}, \cite{Sedghi 2022}, \cite{Ibis 2014}, \cite{Math 2022}.
In this paper, we determined the formulae for calculating $\mathcal{A}$, $V$ and $\mathcal{S}$ in a {\sc H\" ugelsch\" affer} model and we expect that the obtained result will be applicable in various engineering areas.

\bigskip
\noindent
{\bf Acknowledgments.}
The second author was supported in part by the Serbian Ministry of Education, Science and Technological Development,
under project 451-03-68/2022-14/200103.

\newpage

\noindent
{\bf Maja Petrovi\' c}\\
University of Belgrade --\\
The Faculty of Transport and Traffic Engineering,\\
Vojvode Stepe 305,\\
11000 Belgrade, Serbia\\
E-mail: {\it majapet@sf.bg.ac.rs}

\vspace{0.5cc}

\noindent
{\bf Branko Male\v sevi\' c}\\
University of Belgrade,\\
School of Electrical Engineering,\\
Bulevar Kralja Aleksandra 73,\\
11000 Belgrade, Serbia\\
E-mail: {\it branko.malesevic@etf.bg.ac.rs}


\begin{thebibliography}{1}



\bibitem{Kolkata 2013}
{\small {\sc N. B. Basu, A. Dey, D. Ghosh:} {\it Kolkata’s brick sewer renewal: history, challenges and benefits.}
Proceedings of the Institutionof Civil Engineers – Civil Engineering, {\bf 166} (2) (2013), 74--81.}
\bibitem{BEggC 2018}
{\small {\sc J. D. Biggins, J. E. Thompson, T. R. Birkhead:} {\it Accurately quantifying the shape of birds’ eggs.}
Ecology and evolution, {\bf 8} (2018), 9728--9738.}
\bibitem{FlowMI 2011}
{\small {\sc M. Bijankhan, S. Kouchakzadeh:} {\it Egg-shaped cross section: Uniform flow direct solution and stability identification.}
Flow Measurement and Instrumentation, {\bf 22}(6), 2124 (2011), 511--516.}
\bibitem{Boermans 2004}
{\small {\sc R. Berger, L. M. M. Boermans:} {\it Aerodynamic design of the wing-fuselage junction for the high-performance sailplane M\" u-31.}
Technical soaring, {\bf 28}(3) (2004), 13--23.}
\bibitem{Boermans 2006}
{\small {\sc L. M. M. Boermans:} {\it Research on sailplane aerodynamics at Delft University of Technology. Recent and present developments.}
Presented to the Netherlands Association of Aeronautical Engineers (NVvL) on 1 June 2006, 1--~25. \url{http://frotor.fs.cvut.cz/doc/37.pdf}}
\bibitem{Indian 1976}
{\small {\sc Bureau of Indian Standards, IS 4880-2:} {\it  Code of practice for design of tunnels conveying water, Part 2: Geometric design [WRD 14: Water Conductor Systems]} (1976), 1--15.}

\bibitem{Majamath 2019}
{\small {\sc R. Darmawan:} {\it Persamaan Garis Singgung Kurva Bentuk Telur H\" ugelsch\" affer.}
Majamath, {\bf 2}(2) (2019), 94--101.}

\bibitem{XYZ 1944}
{\small {\sc Flugsport:} {\it  Luftfahrt Zeitschrift Flugsport-Kompletter Jahrgang (I. Eierkurven nach F. H\" ugelsch\" affer).}
Digitale Luftfahrt-Bibliothek, {\bf 9} (1944), 162--163.}
\bibitem{Ferreol 2009H}
{\small {\sc R. Ferr\' eol:} {\it Hyperbolism and antihyperbolism of a curve: Newton transformation.} (2009)
\url{https://mathcurve.com/courbes2d/hyperbolisme/hyperbolisme.shtml}
\url{https://mathcurve.com/courbes2d.gb/hyperbolisme/hyperbolisme.shtml}}
\bibitem{Ferreol 2009}
{\small {\sc R. Ferr\' eol:} {\it Oeuf de H\" ugelsch\" affer.} (2009)
\url{https://mathcurve.com/courbes2d/oeuf/oeuf.shtml}}
\bibitem{Ferreol 2001}
{\small {\sc R. Ferr\' eol, J. Mandonnet:} {\it Cubique Mixte.} (2001)
\url{https://mathcurve.com/courbes2d/cubicmixte/cubicmixte.shtml}}


\bibitem{Ryzhik 2015}
{\small {\sc I. Gradshteyn, I. Ryzhik:}{\it Table of Integrals Series and Products.}
8th edn. Academic Press, Cambridge, 2015.}

\bibitem{WaterS 2021}
{\small {\sc R. L. Hachemi, M. Lakehal, B. Achour:} {\it Modern vision for critical flow in an egg-shaped section.}
Water Science \& Technology, {\bf 84}(4), 841 (2021), 840--850.}
\bibitem{Heck 2008}
{\small {\sc A. Heck:} {\it Mathematical Brooding over an Egg.}  Journal of Online Mathematics and Its Applications, Convergence {\bf 8} (2008), 1--12.
\url{https://www.maa.org/external_archive/joma/Volume8/Heck/index.html}}
\bibitem{UTEP 2022}
{\small {\sc S. Holguin, V. Kreinovich:} {\it Shape of an Egg: Towards a Natural Simple Universal Formula.}   Technical Report UTEP-CS-22-55, (2022), 1--6.
\url{https://www.cs.utep.edu/vladik/2022/tr22-55.pdf}}

\bibitem{Cleveland 2008}
{\small {\sc  F. G. Jackson, J. Wasik, R. E. Devaul:} {\it  Standard Construction Drawings. Standard plan for Egg Shaped Sewers Dimensions \& Areas.
City of Cleveland, Ohio, Division of Engineering and Construction}, File No.{\bf 73 M} (2008), 24--24.
\url{https://www.clevelandohio.gov/node/3063};
\url{https://www.clevelandohio.gov/sites/default/files/forms_publications/Standard%20Construction%20Drawings.pdf?id=6510}
 }

\bibitem{VMCh 2021}
{\small {\sc O. Karabulut:} {\it Estimation of the external quality characteristics of goose eggs of known breadth and length.}  Vet Med-Czech, {\bf 66}(10) (2021), 440--447.}

\bibitem{Kostin 1977}
{\small {\sc J. V. Kostin:} {\it On the methods of oomorphological studies and the unifying of the descriptions of oological materials.}
pp. 14-22 in (G. A. Noskov, Ed.): Metodiki issledovanija produktivnosty i struktury vidov ptits v predelah ih arealov, Vilnius, Lithuania, Mokslas 1977.}


\bibitem{EggSurf 2015}
{\small {\sc A. R. Maulana, M. Yunus, D. R. Sulistyaningrum:} {\it The Constructions of Egg-Shaped Surface Equations using H\" ugelsch\" affer's Egg-Shaped Curve.}
Indonesian Journal of Physics, {\bf 26}(2) (2015), 26--30.}
\bibitem{Milojevic 2020}
{\small {\sc M. Milojevi\' c, Z. Joki\' c, S. Mitrovi\' c:} {\it Effects of Morphometric Indicators on Incubation Values of Eggs and Sex of the Chicks of the Light Hen Hybrids.}
In book: Animal Models in Medicine and Biology, {\bf Ch. 13} (2020), 1--15.}

\bibitem{Applet 2022}
{\small {\sc M. Milo\v sevi\' c:} {\em H\" ugelsch\" affer egg curve and surface applet}, {\sc OviForm (2022)} \url{https://oviform.etf.bg.ac.rs} }

\bibitem{EggCurve 2017}
{\small {\sc I. S. Mytiai, A. V. Matsyura:} {\it Geometrical standards in shapes of avian eggs.}
Ukrainian Journal of Ecology, {\bf 7}(3) (2017), 264--282.}
\bibitem{EggCurve 2018}
{\small {\sc I. S. Mytiai, A. V. Matsyura:} {\it Usage of the iterative photo--computing method in specifying of bird egg radiuses curvature.}
Ukrainian Journal of Ecology, {\bf 8}(4) (2018), 195--204.}
\bibitem{EggCurve 2019}
{\small {\sc I. S. Mytiai, A. V. Matsyura:} {\it Mathematical interpretation of artificial ovoids and avian egg shapes (part I).}
Regulatory Mechanisms in Biosystems, {\bf 10}(1) (2019), 92--102.}

\bibitem{Narushin 2020a}
{\small {\sc V. G. Narushin, G. Lu, J. Cugley, M. N. Romanov, D. K. Griffin:} {\it A 2-D imaging-assisted geometrical transformation method for non-destructive evaluation of the volume and surface area of avian eggs.}
Food Control, {\bf 112} (107112) (2020), 1--8.}
\bibitem{Narushin 2020b}
{\small {\sc V. G. Narushin, M. N. Romanov, G. Lu, J. Cugley, D. K. Griffin:} {\it Digital imaging assisted geometry of chicken eggs using H\" ugelsch\" affer's model.}
Biosystems Engineering, {\bf 197} (2020), 45--55.}
\bibitem{Narushin 2021a}
{\small {\sc V. G. Narushin, M. N. Romanov, G. Lu, J. Cugley, D. K. Griffin:} {\it How oviform is the chicken egg? New mathematical insight into the old oomorphological problem.}
Food Control, {\bf 119} (107484) (2021), 1--12.}
\bibitem{Narushin 2021b}
{\small {\sc V. G. Narushin, M. N. Romanov, D. K. Griffin:} {\it Non-destructive measurement of chicken egg characteristics: improved formulae for calculating egg volume and surface area.}
Biosystems Engineering, {\bf 201} (2021), 42--49.}
\bibitem{Narushin 2021c}
{\small {\sc V. G. Narushin, M. N. Romanov, D. K. Griffin:} {\it Non-destructive evaluation of the volumes of egg shell and interior: Theoretical approach.}
Journal of Food Engineering, {\bf 300} (110536) (2021), 1--8.}
\bibitem{Narushin 2021d}
{\small {\sc V. G. Narushin, M. N. Romanov, D. K. Griffin:} {\it Egg and math: introducing a universal formula for egg shape.}
Annals of the New York Academy of Sciences, {\bf 1505} (2021), 169--177.}
\bibitem{Narushin 2022}
{\small {\sc V. G. Narushin, M. N. Romanov, B. Mishra, D. K. Griffin:} {\it Mathematical progression of avian egg shape with associated area and volume determinations.}
Annals of the New York Academy of Sciences, (2022), 1--14.}
\bibitem{Narushin 2022b}
{\small {\sc V. G. Narushin, M. N. Romanov, D. K. Griffin:} {\it  Egg-inspired engineering in the design of thin-walled shelled vessels: a theoretical approach for shell strength.}  Front. Bioeng. Biotechnol. {\bf 10}:995817. (2022),  1--11. doi: 10.3389/fbioe.2022.995817}

\bibitem{Newton 1860}
{\small {\sc Sir I. Newton:}{\it Enumeration of Lines of the third order.}
H. G. Bohn, London, pp. 21--23, 62--63, 1860.}

\bibitem{moNG 2010b}
{\small {\sc M. Obradovi\' c, B. Male\v sevi\' c, M. Petrovi\' c:} {\it Conic Section of a Type of Egg Curve Based Conoid.}
In Proceedings of $2^{nd}$ International Conference for Geometry and Graphics -- moNGeometrija 2010, Belgrade, Serbia (2010), 447--466.}
\bibitem{moNG 2008}
{\small {\sc M. Obradovi\' c, M. Petrovi\' c:} {\it Spatial Interpretation of H\" ugelsch\" affer's Egg Curve Construction.}
In Proceedings of $24^{th}$ National and $1^{st}$ International Conference for Geometry and Graphics -- moNGeometrija 2008, Vrnja\v cka Banja, Serbia (2008), 222--232.}
\bibitem{ICEGD 2013}
{\small {\sc M. Obradovi\' c, B. Male\v sevi\' c, M. Petrovi\' c, G. Djukanovi\' c:} {\it Generating Curves of Higher Order Using the Generalisation of H\" ugelsch\" affer's Egg Curve Construction.}
In Proceedings of International Conference on Engineering Graphics and Design -- ICEGD 2013, Timisoara, Romania, (SCIENTIFIC BULLETIN of the "POLITEHNICA" University of Timisoara, Romania), Tomul 58(72), Fasc. {\bf 1}(2013), 110--115.}

\bibitem{Paganelli 1974}
{\small {\sc C. V. Paganelli, A. J. Olszowka, A. Ar:} {\it  The Avian Egg: Surface Area, Volume, and Density.}
The Condor, {\bf 76}(3)(1974), 319--325. }
\bibitem{MP 2010}
{\small {\sc M. Petrovi\' c: }{\it Egg curves and generalisation H\" ugelsch\" affer's construction.}
M.Sc. degree thesis, Faculty of arhitecture, University of Belgrade, 2010.}
\bibitem{moNG 2010a}
{\small {\sc M. Petrovi\' c, M. Obradovi\' c:} {\it The Complement of the H\" ugelsch\" affer's construction of the Egg Curve.}
In Proceedings of $2^{nd}$ International Conference for Geometry and Graphics -- moNGeometrija 2010, Belgrade, Serbia (2010), 520--531.}
\bibitem{ICEGD 2011}
{\small {\sc M. Petrovi\' c, M. Obradovi\' c, R. Mijailovi\' c:} {\it Suitability analysis of H\" ugelsch\" affer's egg curve application in architectural structures' geometry.}
In Proceedings of International Conference on Engineering Graphics and Design -- ICEGD 2011 "Sustainable Eco Design", Iasi, Romania (Buletinul Institutului Politehnic din Iasi, Publicat de Universitatea Tehnic\v a "Gheorghe Asachi" din Iasi), Tomul LVII(LXI), Fasc. {\bf 3} (2011), 115--122.}

\bibitem{Water 2016}
{\small {\sc M. Regueiro-Picallo, J. Naves, J. Anta, J. Puertas, J. Su\' arez:} {\it Experimental and Numerical Analysis of Egg-Shaped Sewer Pipes Flow Performance.}
Water, {\bf 8}(12), 587 (2016), 1--9.}

\bibitem{Spheroids 1960}
{\small {\sc J. Satterly:} {\it Formulae for Volumes, Surface Areas and Radii of Gyration of Spheres, Ellipsoids and Spheroids.}
The Mathematical Gazette, {\bf 44}(347) (1960), 15--19.}
\bibitem{Schmidbauer 1948}
{\small {\sc H. Schmidbauer: }{\it Kleine Mitteilungen, II. Eine exakte Eierkurvenkonstruktion mit technischen Anwendungen.}
Elemente der Mathematik, {\bf 3} (1948), 67--68.}
\bibitem{Schmidbauer 1949}
{\small {\sc H. Schmidbauer: }{\it Berichtigung.}
Elemente der Mathematik, {\bf 4} (1949), 96--96.}
\bibitem{Sedghi 2022}
{\small {\sc M. Sedghi, M. Ghaderi:} {\it Digital analysis of egg surface area and volume: Effects of longitudinal axis, maximum breadth and weight.}
Accepted in Information Processing in Agriculture, (2022), 1--11.}
\bibitem{Severa 2013}
{\small {\sc L. Severa , \v S. Nedomov\' a , J. Buchar, J. Cupera:} {\it Novel Approaches in Mathematical Description of Hen Egg Geometry.}  International Journal of Food
Properties, {\bf 16}(7) (2013), 1472--1482.}
\bibitem{Water 2019}
{\small {\sc H. Shang, S. Xu, K. Zhang, L. Zhao:} {\it Explicit Solution for Critical Depth in Closed Conduits Flowing Partly Full.}
Water, {\bf 11}(10), 2124 (2019), 1--17.}
\bibitem{WaterSIE 2017}
{\small {\sc N. Stani\' c, J. Langeveld, T. Salet, F. Clemens:} {\it Relating the structural strength of concrete sewer pipes and material properties
retrieved from core samples.}
Structure and Infrastructure Engineering, {\bf 13}(5) (2017), 637--651.}
\bibitem{WaterST 2013}
{\small {\sc N. Stani\' c, C. de Haan, M. Tirion, J.G. Langeveld, F.H.L.R. Clemens:} {\it Comparison of core sampling and visual inspection for assessment of concrete sewer pipe
condition.}
Water Sci. Technol. {\bf 67} (2013), 2458--2466.}

\bibitem{Ibis 2014}
{\small {\sc J. Troscianko:} {\it A Simple Tool for Calculating Egg Shape, Volume and Surface Area from Digital Images.}
Ibis, {\bf 156}(4) (2014), 874--878.}

\bibitem{Tables 2006}
{\small {\sc  HR Wallingford (Firm), D. I. H. Barr:} {\it Tables for the Hydraulic Design of Pipes, Sewers and Channels.}
Thomas Telford, London, Volume {\bf II}, 8th edition (2006), 64--76.}
\bibitem{Heinrich 1908}
{\small {\sc H. Wieleitner: }{\it Spezielle Ebene kurven.}
Leipzig, p. 49, 1908.}
\bibitem{Math 2022}
{\small {\sc Y.-K. Weng, C.-H. Li, C.-C. Lai, C.-W. Cheng:} {\it Equation for Egg Volume Calculation Based on Smart’s Model.}  Mathematics, {\bf 10} 1661 (2022), 1--9.}
\bibitem{WaterER 2022}
{\small {\sc H. Wu, Y. Huang, L. Chen, Y. Zhu, H. Li:} {\it Shape optimization of egg-shaped sewer pipes based on the nondominated sorting genetic algorithm $($NSGA-II$)$.}
Environmental Research, {\bf 204}, 111999 (2022), 1--10.}


\end{thebibliography}
\end{document}